\documentclass[reqno,twoside,final]{amsart}
\usepackage{amsmath}

\newtheorem{Thm}{Theorem}
\newtheorem{Cor}[Thm]{Corollary}

\theoremstyle{definition}

\theoremstyle{remark}
\newtheorem{Rem}{Remark}

\begin{document}
\title[Bounded oscillating linear systems]
{ Systems of linear ordinary differential equations with bounded
coefficients may have very oscillating solutions }

\author{D. Novikov}
\address{Department of Mathematics, Toronto University, Toronto, Canada}
 \email{dmitry@math.toronto.edu}
\abstract
An elementary example shows that the number of zeroes of a component
of a solution of  a system of linear ordinary differential equations
cannot be estimated through the norm of coefficients of the system.
\endabstract
\maketitle

\subsection*{Bounds for oscillation of solutions}
In \cite{IYa} it was shown that a linear ordinary differential
equation of order $n$ with real analytic coefficients bounded in a
neighborhood of the segment $[-1,1]$, admits an upper bound for
the number of isolated zeros for all its solutions.

Actually, the analyticity condition can be relaxed: it is only
the boundedness of the coefficients that matters. Probably, the
simplest result in this spirit looks as follows. Consider the
linear ordinary differential equation
\begin{equation}\label{eqn}
y^{(n)}(t)+a_1(t)y^{(n-1)}(t)+...+a_n(t)y(t)=0.
\end{equation}
with continuous coefficients bounded on
$[\alpha,\beta]\subset\mathbb R$: $|a_i(t)|\le C$ for some $C\ge
1$.

\begin{Thm}[see \cite{NY,Yafields}]
The number of isolated zeros of any solution of the equation
\eqref{eqn} on $[\alpha,\beta]$ cannot exceed $n-1+\frac n {\ln
2} C|\beta-\alpha|$.
\end{Thm}

An analog of this result for \emph{systems} of ordinary
differential equations (vector fields in space) would concern the
number of isolated intersections between integral trajectories of
the field and hyperplanes (and, more generally, hypersurfaces).
Such result was proved in \cite{NY}, see also \cite{lleida}, for
\emph{polynomial} vector fields and algebraic hypersurfaces, and
the bound for the number of isolated intersections is given in
terms of the \emph{height} (maximal magnitude of coefficients) of
the polynomials in the right hand side of the system. Consider a
system of polynomial ordinary differential equations of degree
$d$ in $\mathbb R^n$,
\begin{equation}\label{sys}
  \dot x_i=v_i(t,x),\quad i=1,\dots,n,
  \qquad v_i(t,x)=\sum\nolimits_{k+|\alpha|\le d}v_{ik\alpha}t^kx^\alpha,
\end{equation}
and an algebraic hypersurface $\{P=0\}$, where $P=P(t,x)$ is
another polynomial of the same degree.

\begin{Thm}[see \cite{lleida,NY}]
Assume that all coefficients of the system \eqref{sys} are
bounded by some constant $C$, $|v_{ik\alpha}|\le C$. Let
$\gamma(t)$ be an integral curve of this vector field lying
entirely in the box $B_C=\{|t|<C, |x_i|<C\}\subset\mathbb R^{n+1}$
of the same size $C$.

Then the number of isolated intersections between $\gamma$ and
$\{P=0\}$ can be at most $(2+R)^{B}$, where $B=B(n,d)$ is an
explicit elementary function of $d$ and $n$ only, growing no
faster than $\exp\exp\exp\exp(4n\ln d +O(1))$ as $d,n\to\infty$.
\end{Thm}

As it was remarked in \cite{NY}, this result is nontrivial even
for linear systems
\begin{equation}\label{linsys}
  \dot x=A(t)x,\qquad x\in\mathbb R^n,
  \quad A(t)=\sum\nolimits_{k=0}^d A_k t^k,
\end{equation}
and linear hyperplanes $\{\sum_1^n p_i x_i=0\}$. In this case the
box condition reduces to specifying the interval $t\in[-C,C]$ and
the height condition means that the norms of the matrix
coefficients $A_k\in\operatorname{Mat}_{n\times n}(\mathbb R)$
are assumed to be bounded, $\|A_k\|\le C$.

\begin{Cor}
The number of isolated zeros of any component of a solution for
the linear system \eqref{linsys} on the interval $[-C,C]$ admits
an explicit upper bound in terms of $C$ as above, uniformly over
all solutions of this system.
\end{Cor}

Comparing these two theorems suggests the question whether the
polynomiality condition in the second theorem can be relaxed and
replaced, say, by the norm $\max_{i=1,\dots,n,\,(t,x)\in
B_C}|v_i(t,x)|$ or, in the case of the linear system
\eqref{linsys}, by $\max_{t\in[-C,C]}\|A(t)\|$. We show that
\emph{this is impossible}.

\subsection*{The example}
We construct a linear $2\times 2$-system \eqref{linsys} that is
polynomial of an arbitrarily high degree $2d$, but has the
coefficient matrix $A(t)$ bounded by $1$ in the sense of the norm
on $[-1,1]$ in such a way that the number of isolated zeros of its
first component is $d$ and hence is unbounded in terms of
$\max_{t\in[-C,C]}\|A(t)\|$.

Let $t_1,...,t_d\in[-1,1]$ be any different real numbers from the
interval $[-1,1]$. Consider a linear differential equation of the
first order $\dot x_1=a(t)x_1$, where
$a(t)=\lambda(t-t_1)...(t-t_d)$ and $\lambda$ is chosen so small
that $|a(t)|+|\dot a(t)+a^2(t)|<1$ for any $t\in[-1,1]$.

Its solution $\phi_1=\exp(\int a(t) dt)$ has no zeroes at all, but
the derivative $\phi_2=\dot\phi_1=a(t)\phi_1$ has the same zeroes
as $a(t)$ and satisfies the equation $\dot \phi_2=(\dot
a+a^2)\phi_1$.

Consider now the linear system
$$
\dot x_1=a(t)x_1,\quad\dot x_2=(\dot a(t)+a(t)^2)x_1.
$$

The curve $\gamma=(\phi_1(t), \phi_2(t))$ is a solution of this
system and its second component  has $d$ zeroes on $[-1,1]$. The
coefficients of the system are bounded by $1$ in the supremum
norm on the interval $[-1,1]$. The size of solution is irrelevant
since the system is linear (in other words, we can multiply it by
a constant so small that it will not leave the box $B_1$ for
$|t|\le1$). However, choosing $d$ sufficiently large, we can
obtain any number of zeroes.

\begin{Rem}
The same example also shows that one cannot extend estimates  of
the Theorem~1 to derivatives of the solution.

Actually, by choosing $\lambda$ sufficiently small one can ensure
that the coefficients of the constructed system are uniformly
small in an arbitrarily chosen \emph{complex} neighborhood of the
real segment $[-1,1]$. This shows that the bounds for oscillation
around hyperplanes cannot be achieved in the spirit of \cite{IYa}
(in terms of the bounds for analytic coefficients in the complex
domain) as well.
\end{Rem}

\subsection*{Acknowledgment}
I am grateful to S.~Yakovenko for drawing my attention to this
problem and stimulating discussions.

\bibliographystyle{amsplain}

\begin{thebibliography}{1}
\frenchspacing

\bibitem{IYa}
Yu. Ilyashenko and S. Yakovenko, Counting real zeros of analytic
functions satisfying linear ordinary differential equations,
Journal of Differential equations {\bf 126} (1996), no. 1, 87-105.

\bibitem{lleida}
D. Novikov and S. Yakovenko, Meandering of trajectories of
polynomial vector fields in the affine $n$-space, {\em Publ.
Mat.} {\bf 41} (1997), no.~1, 223--242.

\bibitem{NY}
\bysame, \bysame, Trajectories of polynomial vector fields and
ascending chains of polynomial ideals, Ann. Inst. Fourier {\bf
49} (1999), no.~2, 563--609.

\bibitem{Yafields}
S.~Yakovenko, \emph{On functions and curves defined by ordinary
differential equations}, Proceedings of the Arnoldfest (Ed. by
E.~Bierstone, B.~Khesin, A.~Khovanskii, J.~Marsden), Fields
Institute Communications, 1999, pp.~203--219.

\end{thebibliography}

\end{document}